\documentclass[preprint,12pt]{elsarticle}

\usepackage{amssymb}
\usepackage{units}
\usepackage{amssymb,amsmath}
\usepackage{ntheorem}
\usepackage{graphics,graphicx}
\usepackage{units}
\usepackage{color}
\usepackage{epstopdf}
\usepackage{booktabs}
\usepackage{caption}
\usepackage[justification=centering]{caption}
\usepackage{soul}%高亮头文件
\usepackage{lipsum}
\usepackage{amsfonts}
\usepackage{algorithmic}

\newtheorem{proposition}{Proposition}[section]

\makeatletter

\newcommand{\Rmnum}[1]{\expandafter\@slowromancap\romannumeral #1@}
\makeatother

\journal{Journal}

\begin{document}

\begin{frontmatter}
\title{Block triangular preconditioning for elliptic boundary optimal control with mixed boundary conditions
\footnote{This work was supported by the National Natural Science Foundation of China (No. 12001022).
\\
\indent $^\ast$ Corresponding author: Chaojie Wang.
\\
\indent {\emph{E-mail address: wangcj2019@btbu.edu.cn} }}}
\author{Chaojie Wang$^{a,\ast}$}
\address{a. School of Mathematics and Statistics, Beijing Technology and Business University, Beijing 100048, China}
\begin{abstract}
In this paper, preconditioning the saddle point problem arising from the elliptic boundary optimal control problem with mixed boundary conditions is considered. A block triangular reconditioning method is proposed based on permutations of the saddle point problem and approximations of the corresponding Schur complement.  The spectral properties of the preconditioned matrix is analyzed. Numerical experiments are conducted to demonstrate the effectiveness of the proposed preconditioning method.

\end{abstract}
\end{frontmatter}

\section{Introduction}

The effective solving of optimal control problem with PDE constraints plays an important role in a variety of areas \cite{mgp}-\cite{ftr}. The performance of iterative methods for these problems is normally decided by the choice of preconditioning techniques. In recent years, a lot of preconditioning methods have been presented in various fields \cite{hli}-\cite{cwe}. It has been attracted extensive attentions to propose preconditioning techniques for the structured linear equations appearing in the optimal control problem with different PDE constraints, boundary conditions and control types. In \cite{jsw1}, Schoberl and Zulehner constructed a kind of symmetric and indefinite preconditioners for distributed control problems with elliptic state equations. For the saddle point problems arising from the optimal control problem with a hyperbolic constraint, Benzi, Haber, and Taralli constructed the block triangular preconditioners with diagonal perturbations of the approximate Hessian \cite{mbe2}. As a comparison of preconditioned Krylov subspace iteration methods, Axelsson, Farouq, and Neytcheva constructed a series of block matrix preconditioners for optimal control problems with Possion and convection-diffusion control, respectively \cite{oas1}. Besides, other preconditioning techniques based on domain decomposition and norm equivalence (see, e.g., Heinkenschloss and Nguyen \cite{mhh1},  Arioli, Kourounis, and Loghin \cite{mad2}), operator methods (see, e.g., Zulehner \cite{wz1}, Elvetun and Nielsen \cite{oeb1}) have also been investigated on solving the PDE-constrained optimization problem more efficiently.

In this paper, we consider building block preconditioners for the elliptic boundary optimal control problem
\begin{equation}\label{eq1}
\min_{y,u}\frac{1}{2}||y-y_{d}||^{2}_{L_{2}(\Omega)}+\frac{\beta}{2}||u||^{2}_{L_{2}(\Gamma)},
\end{equation}
subject to
\begin{equation}\label{eq2}
\renewcommand\arraystretch{1.5}
\left\{
\begin{array}{cl}
-\nabla^{2}y=f & \text{in}\quad\Omega,\\
y=0                   & \text{on}\quad\Gamma^{c},\\
\partial_{n} y=u & \text{on}\quad\Gamma,\\
\end{array}
\right.
\end{equation}
where $\Omega\subset\mathbb{R}^{2}$, $\overline{\Gamma^{c}}\cup \overline{\Gamma}=\overline{\partial\Omega}, \Gamma^{c}\cap\Gamma=\emptyset$, $y$ denotes the state variable, $u$ is the control variable and $\beta$ denotes the positive regularization parameter. The solution for $y$ ought to approach $y_{d}$ as near as possible under the listed PDE constraints. 

Over the last decade, several preconditioning methods have been proposed for the saddle point problem arising from this problem. For the pure Neumann boundary case ($\Gamma^{c}=\emptyset$), Rees, Dollar and Wathen just used one term of the Schur complement as its approximation and used the strategy to build block preconditioners \cite{trh1}. The resulting iteration number varied as the mesh size decreased. In \cite{jpm}, Pearson and Wathen built preconditioners for the same problem based on a matching strategy, which was originally proposed by them for constructing preconditioners for the distributed Poisson control problem \cite{jpa}. In the distributed Poisson problem, the matching strategy renders the eigenvalues of the preconditioned matrix  bounded by certain intervals independent of the mesh size and the regularization parameter. However, the eigenvalue properties are no longer valid for the optimal control with pure Neumann boundaries and the iteration number was found to be mesh size dependent. It is noticed that little effort has been directed at preconditioning the saddle point problem arising from the mixed boundary case ($\Gamma^{c}\neq\emptyset$). This will be the focus in this paper. Our preconditioner will also take advantage of the structure of the saddle point problem; however, we choose to work with a permuted version which yields a different sparse approximation for the Schur complement. The resulting preconditioner performs well for a practical range of the regularization parameter and exhibits better mesh independence.

The rest of this paper is organized as the following. In Section 2, we derive the saddle point problem corresponding to the elliptic boundary optimal control problem with mixed boundary conditions and give a brief description of related work on this issue. In Section 3, we design a block triangular preconditioner and give the analysis of certain spectral properties for the preconditioned system. And numerical results provided in Section 4 illustrate effectiveness of our proposed preconditioning method. Finally, we give some concluding remarks as the last section.

\section{Background}
The weak form of (\ref{eq1})-(\ref{eq2}) reads : find $(y, u)\in H^{1}_{\Gamma^{c}}(\Omega)\times L_{2}(\Gamma)$ to the problem
\begin{equation}\label{eq3}
\renewcommand\arraystretch{2}
\left\{
\begin{array}{l}
\displaystyle\min\frac{1}{2}||y-y_{d}||^{2}_{L_{2}(\Omega)}+\frac{\beta}{2}||u||^{2}_{L_{2}(\Gamma)},\\
\displaystyle\text{s.t.}\int_{\Omega}\nabla y\cdot\nabla z d\Omega-\int_{\Gamma}u z ds=\int_{\Omega}f z d\Omega, \quad \forall z\in H^{1}_{\Gamma^{c}}(\Omega),\\
\end{array}
\right.
\end{equation}
where $H^{1}_{\Gamma^{c}}(\Omega)=\{z: z\in H^{1}(\Omega), z|_{\Gamma^{c}}=0\}$. Normally, both the \emph{discretize-then-optimize} strategy and \emph{optimize-then-discretize} strategy can be used to tackle PDE-constrained  optimal control problems. As for the optimal control problem constrained by the Poisson equation, these two strategies result in the same linear system \cite{yq1}. Here we choose the \emph{discretize-then-optimize} strategy combined with the finite element method and the Lagrange multiplier method to solve the elliptic boundary control problem.

\subsection{Discretization}
Let $V^{h}\subset H^{1}_{\Gamma^{c}}(\Omega)$ be a trail space and its basis is $\{\phi_{1}, \cdots,\phi_{n}\}$. We assume that $n_{I}$ basis elements have nonzero support inside $\Omega$ with the remaining $n_{B}$ basis elements having support on $\Gamma$ and $n=n_{I}+n_{B}$. Then any $y_{h}\in V^{h}$ has the form
\begin{equation}\label{eq4}
y_{h}=\sum^{n}_{i=1}Y_{i}\phi_{i},
\end{equation}
where $\bold{y}=(Y_{1},\cdots,Y_{n})^{T}$ denotes the coefficient vector.
Let $W^{h}\subset L_{2}(\Gamma)$ be a test space and its basis is $\{\psi_{1},\cdots,\psi_{m_{B}}\}$ with nonzero support on the Neumann boundary $\Gamma$. Then any $u_{h}\in W^{h}$ has the form
\begin{equation}\label{eq5}
u_{h}=\sum^{m_{B}}_{j=1}U_{j}\psi_{j},
\end{equation}
where $\bold{u}=(U_{1},\cdots,U_{m})^{T}$ denotes the coefficient vector . 

According to the Galerkin finite element method, the problem (\ref{eq3}) can be discretized as:  find $(y_{h},u_{h})\in V_{h}\times W_{h}$ to the problem
\begin{equation}\label{eq6}
\renewcommand\arraystretch{2}
\left\{
\begin{array}{l}
\displaystyle\min\frac{1}{2}||y_{h}-y_{d}||^{2}_{L_{2}(\Omega)}+\frac{\beta}{2}||u_{h}||^{2}_{L_{2}(\Gamma)},\\
\displaystyle\text{s.t.}\int_{\Omega}\nabla y_{h}\cdot\nabla z_{h} d\Omega-\int_{\Gamma}u_{h} z_{h} ds=\int_{\Omega}f z_{h} d\Omega, \quad \forall z_{h}\in V^{h}(\Omega).\\
\end{array}
\right.
\end{equation}
Substituting (\ref{eq4}) and (\ref{eq5}) into (\ref{eq6}) results in the matrix form
\begin{equation}\label{eq7}
\renewcommand\arraystretch{2}
\left\{
\begin{array}{l}
\displaystyle\min_{\bold{y},\bold{u}}\frac{1}{2}\bold{y}^{T}M\bold{y}-\bold{y}^{T}\bold{b}+
\frac{\beta}{2}\bold{u}^{T}M_{\Gamma}\bold{u},\\
\displaystyle\text{s.t.}K\bold{y}-N_{\Gamma}\bold{u}=\bold{f}, \\
\end{array}
\right.
\end{equation}
where the elements of the mass matrix $M$ and the stiffness matrix $K$ are
$$M_{ij}=\int_{\Omega}\phi_{i}\phi_{j} d\Omega, \quad K_{ij}=\int_{\Omega}\nabla\phi_{i}\cdot\nabla\phi_{j} d\Omega, i, j=1,\cdots,n,$$
the elements of the boundary mass matrix $M_{\Gamma}$ and the matrix $N_{\Gamma}$ are
$$[M_{\Gamma}]_{ij}=\int_{\Gamma}\psi_{i}\psi_{j} d\Omega, \quad [N_{\Gamma}]_{kj}=\int_{\Gamma}\phi_{k}\psi_{j} ds, i, j=1,\cdots,m_{B}, k=1,\cdots,n, $$
and the elements of the vectors $\bold{b}$ and $\bold{f}$ are
$$b_{j}=\int_{\Omega}y_{d}\phi_{j} d\Omega,\quad f_{j}=\int_{\Omega}f\phi_{j} d\Omega, j=1,\cdots,n.$$

Note that the Lagrangian function of (\ref{eq3}) takes the form
$$
L(\bold{y},\bold{u},\bold{p})=\frac{1}{2}\bold{y}^{T}M\bold{y}-\bold{y}^{T}\bold{b}+
\frac{\beta}{2}\bold{u}^{T}M_{\Gamma}\bold{u}+\bold{p}^{T}(K\bold{y}-N_{\Gamma}\bold{u}-\bold{f}),
$$
where $\bold{p}=(P_{1},\cdots,P_{n})$ and $p_{h}=\sum^{n}_{j=1}P_{j}\phi_{j}$ is the finite element approximation of the Lagrange multiplier $p$. The first order optimality conditions give the linear system 
\begin{equation}\label{eq8}
\left(
\begin{array}{cc|c}
M & 0           & K\\
0 & \beta M_{\Gamma} & -N^{T}_{\Gamma}\\
\hline
K & -N_{\Gamma}      & 0
\end{array}
\right)
\left(
\begin{array}{c}
\bold{y}\\
\bold{u}\\
\bold{p}
\end{array}
\right)
=
\left(
\begin{array}{c}
\bold{b}\\
\bold{0}\\
\bold{f}
\end{array}
\right).
\end{equation}

It is noted that in the pure Neumann boundary case $\Gamma^{c}=\emptyset$ and the mixed boundary case $\Gamma^{c}\neq\emptyset$, the corresponding linear systems both posses the saddle point structure shown in (\ref{eq8}). However, the block matrices $M_{\Gamma}$ and $K$ are different. In this paper, we focus on building preconditioners for the saddle point problem arising from the mixed boundary case. 

\subsection{Previous work}
It is noticed that little work has been done to construct preconditioners for the saddle point problem (\ref{eq8}) related to the mixed boundary case. Instead, most effort has been directed at the pure Neumann boundary case, for which block preconditioners were designed based on the 2-by-2 block structure 
$$
\mathcal{A}=\left(
\begin{array}{cc}
A & B\\
B^{T} & 0
\end{array}
\right)
$$
with
$$A=\left(
\begin{array}{cc}
M & 0\\
0 & \beta M_{\Gamma}
\end{array}
\right),
B=\left(
\begin{array}{c}
K\\
-N^{T}_{\Gamma}
\end{array}
\right).
$$
The corresponding Schur complement is 
$$S=KM^{-1}K+\frac{1}{\beta}N_{\Gamma}M^{-1}_{\Gamma}N^{T}_{\Gamma}.$$
In this case, the stiffness matrix $K$ is singular. Here we give a brief description of the work done by Rees \cite{trh1} and Pearson \cite{jpm} for this case.  The resulting preconditioners have block diagonal form and are denoted by
$$
\mathcal{P}_{D,i}=\left(
\begin{array}{cc}
A & 0\\
0  & S_{i}
\end{array}
\right),
i=1,2,
$$
where the approximations $S_{i}$ to the Schur complement are described below. In \cite{trh1},  the Schur complement $S$ was approximated  as
$$S_{1}=KM^{-1}K.$$
In \cite{jpm}, the Schur complement $S$ was approximated  as
$$S_{2}=\bigg(K+\sqrt{\frac{h}{\beta}}N_{\Gamma}M^{-1}_{\Gamma}N^{T}_{\Gamma}\bigg)M^{-1}\bigg(K+\sqrt{\frac{h}{\beta}}N_{\Gamma}M^{-1}_{\Gamma}N^{T}_{\Gamma}\bigg),$$
where $h$ denotes the mesh size. This is the so-called matching strategy and was originally designed for distributed control, where the eigenvalues of the preconditioned Schur complement were shown to be bounded independently of parameters \cite{jpa}. However, for the problem under consideration, no precise spectral bounds hold. Indeed, experiments in \cite{jpm} show that the performance of MINRES deteriorates under mesh refinement and also with reducing regularization parameter.

\section{Preconditioning}
In this section, we design a block triangular preconditioner for the saddle point problem in the elliptic boundary control problem with mixed boundary conditions. And the spectral property of the preconditioned matrices is analyzed.

\subsection{Block triangular preconditioner}
The design of the preconditioner will take advantage of the block structure given in (\ref{eq8}). First, we permute the original saddle point problem into its equivalent form
\begin{equation}\label{eq12}
\renewcommand\arraystretch{1.2}
\underbrace{
\left(
\begin{array}{cc|c}
K & -N_{\Gamma}     & 0\\
0 & \beta M_{\Gamma} & -N^{T}_{\Gamma}\\
\hline
M & 0           & K
\end{array}
\right)
}_{\mathcal{A}}
\left(
\begin{array}{c}
\bold{y}\\
\bold{u}\\
\bold{p}
\end{array}
\right)
=
\left(
\begin{array}{c}
\bold{f}\\
\bold{0}\\
\bold{b}
\end{array}
\right).
\end{equation}
The coefficient matrix of this linear system has the block structure
$$
\mathcal{A}=\left(
\begin{array}{cc}
A & B\\
C & D
\end{array}
\right)
$$
with
$$A=\left(
\begin{array}{cc}
K & -N_{\Gamma}\\
0 & \beta M_{\Gamma}
\end{array}
\right),
B=\left(
\begin{array}{c}
0\\
-N^{T}_{\Gamma}
\end{array}
\right),
C=\left(
\begin{array}{cc}
M & 0
\end{array}
\right),
D=K,
$$
where $A\in \mathbb{R}^{(n+m_{B})\times (n+m_{B})}$, $B\in \mathbb{R}^{(n+m_{B})\times n}$, $C\in \mathbb{R}^{n\times (n+m_{B})}$, $D\in \mathbb{R}^{n\times n}$.
The resulting Schur complement of the linear system takes the form
\begin{equation}\label{eq13}
S=K+\frac{1}{\beta}MK^{-1}N_{\Gamma}M^{-1}_{\Gamma}N^{T}_{\Gamma}.
\end{equation}

Here the Schur complement $S$ is approximated as
\begin{equation}\label{eq14}
\widehat{S}=K.
\end{equation}
Then we can construct a block triangular preconditioner for $\mathcal{A}$ as
\begin{equation}\label{eq15}
\renewcommand\arraystretch{1.2}
\mathcal{P}_{T}=\left(
\begin{array}{cc}
A & B\\
0 & \widehat{S}
\end{array}
\right)
=\left(
\begin{array}{cc|c}
K & -N_{\Gamma}      & 0\\
0 & \beta M_{\Gamma} & -N^{T}_{\Gamma}\\
\hline
0 & 0           & K
\end{array}
\right).
\end{equation}
The linear system $\mathcal{P}_{T}\bold{v}=\bold{d}$ with $\bold{v}=(\bold{v}^{T}_{1}, \bold{v}^{T}_{2}, \bold{v}^{T}_{3})^{T}$ and $\bold{d}=(\bold{d}^{T}_{1}, \bold{d}^{T}_{2}, \bold{d}^{T}_{3})^{T}$ can be solved step by step as follows
\begin{equation}\label{eq16}
\left\{
\begin{array}{c}
\begin{array}{lll}
K\bold{v}_{3}&=&\bold{d}_{3},\\
\beta M_{\Gamma}\bold{v}_{2}&=&\bold{d}_{2}+N^{T}_{\Gamma}\bold{v}_{3},\\
K\bold{v}_{1}&=&\bold{d}_{1}+N_{\Gamma}\bold{v}_{2}.\\
\end{array}
\end{array}
\right.
\end{equation}
As the stiffness matrix $K$ and the boundary mass matrix $M_{\Gamma}$ in the mixed boundary case are nonsingular and sparse, the solutions of the linear systems in (\ref{eq16}) can be obtained at low cost using direct methods or iterative methods.

\subsection{Analysis of the preconditioner}
We now analyze the spectrum of the preconditioned matrix $\mathcal{P}^{-1}_{T}\mathcal{A}$, which is the same as that of $\mathcal{A}\mathcal{P}^{-1}_{T}$. Note that
$$
\mathcal{A}\mathcal{P}^{-1}_{T}=\left(
\begin{array}{cc}
I              & 0\\
CA^{-1} &S\widehat{S}^{-1}
\end{array}
\right),
$$
thus $\mu=1$ is an eigenvalue with multiplicity $n+m_{B}$ and with eigenvectors $e_{j}\in \mathbb{R}^{2n+m_{B}}, j=1,\cdots, n+m_{B}$. 

As for the other eigenvalues, we consider the eigenvalue problem
\begin{equation}\label{eq17}
S\bold{x}=\mu \widehat{S}\bold{x}, \bold{x}\neq 0.
\end{equation}
Denote $K_{M}=K^{-1}MK^{-1}$ and $M_{\gamma}=N_{\Gamma}M^{-1}_{\Gamma}N^{T}_{\Gamma}$, then we have
$$(I+\frac{1}{\beta}K_{M}M_{\gamma})\bold{x}=\mu\bold{x}.$$
If $\bold{x}\in \text{ker}(M_{\gamma})$, then
$$\bold{x}=\mu\bold{x}.$$
Thus $\mu=1$ is an eigenvalue of $S\widehat{S}^{-1}$ with multiplicity $n-m_{B}$ and with eigenvectors $\bold{x}\in \text{ker}(M_{\gamma})$. Combined with above result, we observe that $1$ is an eigenvalue of the preconditioned matrix $\mathcal{P}^{-1}_{T}\mathcal{A}$ with multiplicity $2n$.

Let $\bold{x}\notin \text{ker}(M_{\gamma})$ denote an eigenvector corresponding to the eigenvalue $\mu=\mu_{0}$, we have
$$(I+\frac{1}{\beta}K_{M}M_{\gamma})\bold{x}=\mu_{0} \bold{x}.$$
Then
$$(K^{-1}_{M}+\frac{1}{\beta}M_{\gamma})\bold{x}=\mu_{0} K^{-1}_{M}\bold{x}.$$
Hence,
\begin{equation}\label{eq18}
\mu_{0}=1+\frac{\bold{x}^{T}M_{\gamma}\bold{x}}{\beta \bold{x}^{T}K^{-1}_{M}\bold{x}}.
\end{equation}
 Notice that if the nodes are ordered in a way that all the interior nodes are followed by the boundary nodes, then $M_{\gamma}$ has the block structure
$$
M_{\gamma}=
\left(
\begin{array}{cc}
\bold{0}&\bold{0}\\
\bold{0}&M_{\Gamma}
\end{array}
\right).
$$
Denote $\bold{x}=(\bold{0}, \boldsymbol{z}^{T})^{T}, \bold{\boldsymbol{z}}\neq0$, then
$$\frac{\bold{x}^{T}M_{\gamma}\bold{x}}{\bold{x}^{T}\bold{x}}=\frac{\bold{\boldsymbol{z}}^{T}M_{\Gamma}\bold{\boldsymbol{z}}}{\bold{\boldsymbol{z}}^{T}\bold{\boldsymbol{z}}}.$$
The boundary mass matrix is normally considered in lumped form, which approximates the scaled identity matrix in the sense that \cite{hed1}
\begin{equation}\label{eq19}
c_{1}h\leq\frac{\bold{z}^{T}M_{\Gamma}\bold{z}}{\bold{z}^{T}\bold{z}}\leq c_{2}h, \forall \bold{z}\neq 0\in \mathbb{R}^{m_{B}},
\end{equation}
where $c_{1}$ and $c_{2}$ are constants independent of the mesh size $h$. Besides,
$$\frac{\bold{x}^{T}K^{-1}_{M}\bold{x}}{\bold{x}^{T}\bold{x}}=\frac{\bold{x}^{T}KM^{-1}K\bold{x}}{\bold{x}^{T}\bold{x}}=\frac{(K\bold{x})^{T}M^{-1}(K\bold{x})}{(K\bold{x})^{T}K^{-2}(K\bold{x})}=\frac{\boldsymbol{r}^{T}M^{-1}\boldsymbol{r}}{\boldsymbol{r}^{T}\boldsymbol{r}}\centerdot\frac{\boldsymbol{r}^{T}\boldsymbol{r}}{\boldsymbol{r}^{T}K^{-2}\boldsymbol{r}},$$
where $\boldsymbol{r}=K\bold{x}\neq 0$. Recall that the mass matrix $M$ and stiffness matrix $K$ have the following spectral properties \cite{hed1}
$$c_{3}h^{2}\leq\frac{\bold{x}^{T}M\bold{x}}{\bold{x}^{T}\bold{x}}\leq c_{4}h^{2}, \forall \bold{x}\neq 0\in \mathbb{R}^{n},$$
and
$$d_{1}h^{2}\leq\frac{\bold{x}^{T}K\bold{x}}{\bold{x}^{T}\bold{x}}\leq d_{2}, \forall \bold{x}\neq 0\in \mathbb{R}^{n},$$
where $c_{3}, c_{4}, d_{1}, d_{2}$ are constants independent of the mesh size $h$. Therefore,
\begin{equation}\label{eq20}
\frac{d^{2}_{1}}{c_{4}}h^{2}\leq\frac{\bold{x}^{T}K^{-1}_{M}\bold{x}}{\bold{x}^{T}\bold{x}}\leq\frac{d^{2}_{2}}{c_{3}}h^{-2}.
\end{equation}
Combining (\ref{eq18}), (\ref{eq19}) and (\ref{eq20}), we have
\begin{equation}\label{eq21}
1+\frac{1}{\beta}ch^{3}\leq\mu_{0}\leq 1+ \frac{1}{\beta}dh^{-1},
\end{equation}
where $c=\frac{c_{1}c_{3}}{d^{2}_{2}}$ and $d=\frac{c_{2}c_{4}}{d^{2}_{1}}$ are constants independent of $\beta$ and $h$. 

The above results on eigenvalue properties of the preconditioned matrix $\mathcal{P}^{-1}_{T}\mathcal{A}$ are summarized in Proposition 3.1.
\begin{proposition}\label{prop3-1}
If the block triangular matrix $\mathcal{P}_{T}$ in (\ref{eq15}) is taken as a preconditioner for the matrix $\mathcal{A}$ in (\ref{eq12}), then 1 is an eigenvalue of the preconditioned matrix $\mathcal{P}^{-1}_{T}\mathcal{A}$ with $2n$ multiplicity and the other $m_{B}$ eigenvalues are bounded by $[1+\frac{1}{\beta}ch^{3}, 1+\frac{1}{\beta}dh^{-1}]$, where $c$ and $d$ are constants independent of $\beta$ and $h$.
\end{proposition}

According to Proposition \ref{prop3-1}, the eigenvalue interval grows as $\beta\rightarrow 0$ and $h\rightarrow 0$.  Notice that the preconditioned linear system is non-symmetric and the GMRES method is used in this paper. Unlike the MINRES method, the spectral analysis is not sufficient to determine the convergence of the GMRES method \cite{agv}-\cite{mbg}. The eigenvalue bound provided here, while, can help us to gain a better insight into the property of the proposed preconditioner. Actually, the numerical results in the next section show the good performance of the proposed preconditioner.

\section{Numerical results}
Numerical experiments are conducted in this section to show the effectiveness of the proposed preconditioning method. We solve the elliptic boundary control problem
$$
\min_{y,u}\frac{1}{2}||y-y_{d}||^{2}_{L_{2}(\Omega)}+\frac{\beta}{2}||u||^{2}_{L_{2}(\Gamma)},
$$
subject to
$$
\renewcommand\arraystretch{1.5}
\left\{
\begin{array}{cl}
-\nabla^{2}y=f & \text{in}\quad\Omega,\\
y=0                   & \text{on}\quad\Gamma^{c},\\
\partial_{n} y=u & \text{on}\quad\Gamma.\\
\end{array}
\right.
$$
In all examples, the domain was $\Omega=[0,1]^{2}$ and the forcing term was $f=0$. The regularization parameter was taken to be $\beta=10^{-2}, 10^{-4}, 10^{-6}$ and the boundary $\Gamma$ was taken to be one of the following sets
$$\Gamma_{1}=\{(x, y)\in\partial\Omega: y=1\},$$
$$\Gamma_{2}=\{(x, y)\in\partial\Omega: x=1\}\cup\Gamma_{1},$$
$$\Gamma_{3}=\partial\Omega\backslash\{(x, y)\in\partial\Omega: x<\frac{1}{2}, y<\frac{1}{2}\}.$$
We tested the proposed triangular preconditioner
$$
\renewcommand\arraystretch{1.2}
\mathcal{P}_{T}=\left(
\begin{array}{cc|c}
K & -N_{\Gamma}      & 0\\
0 & \beta M_{\Gamma} & -N^{T}_{\Gamma}\\
\hline
0 & 0           & K
\end{array}
\right)
$$
as a preconditioner for the GMRES method with tolerance of $10^{-6}$ for the relative residual measured in the 2-norm. We compared the performance with that of the MINRES method equipped with the block diagonal preconditioners
$$
\renewcommand\arraystretch{1.2}
\mathcal{P}_{D,i}
=\left(
\begin{array}{cc|c}
M & 0           & 0\\
0 & \beta M_{\Gamma} & 0\\
\hline
0 & 0      & S_{i}
\end{array}
\right),
$$
where
$$S_{1}=KM^{-1}K,$$
$$S_{2}=\bigg(K+\sqrt{\frac{h}{\beta}}N_{\Gamma}M^{-1}_{\Gamma}N^{T}_{\Gamma}\bigg)M^{-1}\bigg(K+\sqrt{\frac{h}{\beta}}N_{\Gamma}M^{-1}_{\Gamma}N^{T}_{\Gamma}\bigg).$$
The results presented here were obtained using a P1P1 discretization. Other discretizations (P1P1, P2P0, P1P0) were also tested with little difference observed in performance. In the following result tables, the iteration numbers and the CPU time (in seconds) in brackets are listed. All numerical tests were conducted using MATLAB with R2016b version on a computer with Intel Core i5-4590 at 3.30 GHz and 8 GB of RAM.

\subsection{Test problem 1}
Consider the desired state
$$
\renewcommand\arraystretch{1.5}
y_{d}=\left\{
\begin{array}{c}
\begin{array}{ccl}
1,&x\leq\frac{1}{2}, y\leq\frac{1}{2},\\
0,&\text{elsewhere}.
\end{array}
\end{array}
\right.
$$
An illustration of the computed state and control corresponding to $\beta=10^{-6}$ and $\Gamma=\Gamma_{3}$ is shown in Fig. 1. The performance comparisons of the preconditioners are presented in Table 1-3.
%figure1
\begin{figure}[t!]
\centering
\centerline{\includegraphics{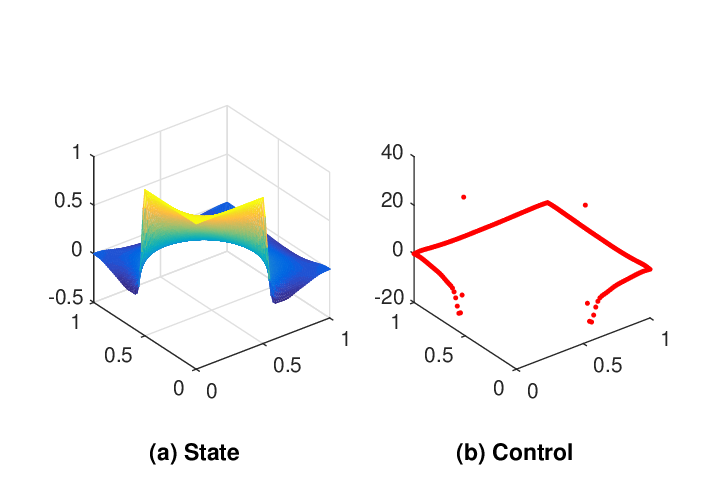}}
\caption*{Fig. 1. Computed state (a) and control (b) for Test problem 1 with $\beta=10^{-6}$ and $\Gamma=\Gamma_{3}$.}
\end{figure}

%table1
\begin{table}[!b]
\renewcommand\arraystretch{1.2}
\scriptsize
\centering
\caption{\protect\\Iterative results for Test problem 1 with $\Gamma=\Gamma_{1}$}
\begin{tabular}{c|ccc|ccc|ccc}
\hline
Method          & \multicolumn{3}{c|}{MINRES($\mathcal{P}_{D,1}$) \cite{trh1}} &  \multicolumn{3}{c|}{MINRES($\mathcal{P}_{D,2}$) \cite{jpm}} &  \multicolumn{3}{c}{$\bold{GMRES(\mathcal{P}_{T}}$)}     \\
\hline
DoF/$\beta$& $10^{-2}$  & $10^{-4}$               & $10^{-6}$  & $10^{-2}$  & $10^{-4}$               & $10^{-6}$  & $10^{-2}$  & $10^{-4}$& $10^{-6}$\\
\hline
4257    &9(0.32)&21(0.40)&68(1.35)&23(1.15)&45(2.10)&63(2.92)&$\bold{3(0.14)}$&	$\bold{5(0.07)}$&	$\bold{7(0.09)}$\\
16705  &9(1.18)&19(1.98)&59(6.06)&25(5.17)&49(10.4)&73(14.2)&$\bold{3(0.32)}$&	$\bold{5(0.48)}$&	$\bold{7(0.53)}$\\
66177  &9(4.57)&15(7.54)&53(25.2)&27(23.3)&57(48.8)&85(72.0)&$\bold{3(1.15)}$&	$\bold{5(1.71)}$&	$\bold{6(2.07)}$\\
263425&9(20.6)&15(33.6)&47(100)&27(113)&57(237)&88(368)&$\bold{3(4.89)}$&$\bold{4(6.14)}$&	$\bold{6(9.06)}$\\
\hline
\end{tabular}
\end{table}

%table2
\begin{table}[!t]
\renewcommand\arraystretch{1.2}
\scriptsize
\centering
\caption{\protect\\Iterative results for Test problem 1 with $\Gamma=\Gamma_{2}$}
\begin{tabular}{c|ccc|ccc|ccc}
\hline
Method          & \multicolumn{3}{c|}{MINRES($\mathcal{P}_{D,1}$) \cite{trh1}} &  \multicolumn{3}{c|}{MINRES($\mathcal{P}_{D,2}$) \cite{jpm}} & \multicolumn{3}{c}{$\bold{GMRES(\mathcal{P}_{T}}$)}     \\
\hline
DoF/$\beta$& $10^{-2}$  & $10^{-4}$               & $10^{-6}$  & $10^{-2}$  & $10^{-4}$               & $10^{-6}$  & $10^{-2}$  & $10^{-4}$& $10^{-6}$\\
\hline
4289    &11(0.21)&28(0.56)&114(2.25)&29(1.71)&55(2.90)&75(3.85)&$\bold{4(0.10)}$&	$\bold{6(0.10)}$&	$\bold{9(0.14)}$\\
16769  &11(1.25)&29(3.06)&101(10.3)&29(7.29)&66(16.1)&89(21.6)&$\bold{4(0.31)}$&	$\bold{6(0.46)}$&	$\bold{8(0.65)}$\\
66305  &11(5.75)&23(11.3)&97(46.6)&33(32.5)&76(73.1)&105(101)&$\bold{4(1.43)}$&	$\bold{5(1.67)}$&	$\bold{7(2.28)}$\\
263681&11(25.4)&23(50.2)&83(175)&35(174)&87(434)&125(608)&$\bold{4(6.28)}$&$\bold{	5(7.84)}$&	$\bold{7(9.87)}$\\
\hline
\end{tabular}
\end{table}

%table3
\begin{table}[!t]
\renewcommand\arraystretch{1.2}
\scriptsize
\centering
\caption{\protect\\Iterative results for Test problem 1 with $\Gamma=\Gamma_{3}$}
\begin{tabular}{c|ccc|ccc|ccc}
\hline
Method          & \multicolumn{3}{c|}{MINRES($\mathcal{P}_{D,1}$) \cite{trh1}} &  \multicolumn{3}{c|}{MINRES($\mathcal{P}_{D,2}$) \cite{jpm}} &  \multicolumn{3}{c}{$\bold{GMRES(\mathcal{P}_{T}}$)}     \\
\hline
DoF/$\beta$& $10^{-2}$  & $10^{-4}$               & $10^{-6}$  & $10^{-2}$  & $10^{-4}$               & $10^{-6}$  & $10^{-2}$  & $10^{-4}$& $10^{-6}$\\
\hline
4321    &19(0.42)&55(1.79)&299(5.82)&47(2.89)&87(5.01)&127(7.40)&$\bold{5(0.14)}$&	$\bold{9(0.12)}$&	$\bold{25(0.37)}$\\
16833  &17(1.81)&55(5.90)&365(36.7)&49(14.0)&115(33.0)&185(51.8)&$\bold{5(0.40)}$&	$\bold{9(0.64)}$&	$\bold{22(1.62)}$ \\
66433  &17(8.48)&57(27.7)&393(184)&53(58.4)&139(151)&257(278)&$\bold{5(1.78)}$&	$\bold{8(2.62)}$&	$\bold{17(5.25)}$\\
263937&15(33.5)&55(115)&393(820)&54(360)&163(1074)&348(2332)&$\bold{4(6.12)}$&	$\bold{8(11.0)}$&	$\bold{14(18.5)}$\\
\hline
\end{tabular}
\end{table}

It can be seen from these tables that the GMRES($\mathcal{P}_{T}$) method requires much less iterations and CPU time than the MINRES($\mathcal{P}_{D,1}$) and MINRES($\mathcal{P}_{D,2}$)  method. The iteration number of the GMRES($\mathcal{P}_{T}$) method exhibits independence on the mesh size, while this is not the case for the other methods. Besides, the required iteration number of the GMRES($\mathcal{P}_{T}$) method keeps reasonable even though all the methods show dependence on the regularization parameter $\beta$. For all the methods, the CPU time increases with the problem size nearly in a linear way. Moreover, it is observed that the iteration number of all the methods grows as the Neumann boundary size becomes larger.

\subsection{Test problem 2}
Consider the desired state
$$
\renewcommand\arraystretch{1.5}
y_{d}=\left\{
\begin{array}{c}
\begin{array}{ccl}
(2x-1)^{2}(2y-1)^2,&x\leq\frac{1}{2}, y\leq\frac{1}{2},\\
0,&\text{elsewhere}.
\end{array}
\end{array}
\right.
$$
An illustration of the computed state and control corresponding to $\beta=10^{-6}$ and $\Gamma=\Gamma_{3}$ is shown in Fig. 2. 
%figure1
\begin{figure}[b!]
\centering
\centerline{\includegraphics{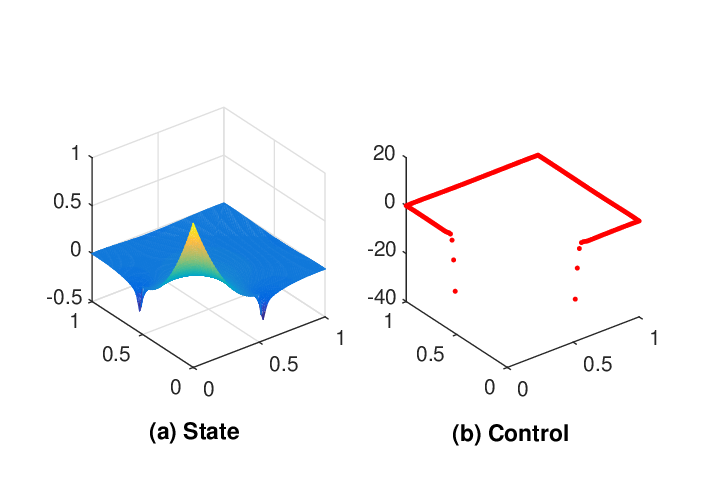}}
\caption*{Fig. 2. Computed state (a) and control (b) for Test problem 2 with $\beta=10^{-6}$ and $\Gamma=\Gamma_{3}$.}
\end{figure}
%table4
\begin{table}[!t]
\renewcommand\arraystretch{1.2}
\scriptsize
\centering
\caption{\protect\\Iterative results for Test problem 2 with $\Gamma=\Gamma_{1}$}
\begin{tabular}{c|ccc|ccc|ccc}
\hline
Method          & \multicolumn{3}{c|}{MINRES($\mathcal{P}_{D,1}$) \cite{trh1}} &  \multicolumn{3}{c|}{MINRES($\mathcal{P}_{D,2}$) \cite{jpm}} &  \multicolumn{3}{c}{$\bold{GMRES(\mathcal{P}_{T}}$)}     \\
\hline
DoF/$\beta$& $10^{-2}$  & $10^{-4}$               & $10^{-6}$  & $10^{-2}$  & $10^{-4}$               & $10^{-6}$  & $10^{-2}$  & $10^{-4}$& $10^{-6}$\\
\hline
4257    &9(0.23)&19(0.43)&63(1.21)&21(1.04)&41(1.92)&55(2.60)&$\bold{4(0.08)}$&	$\bold{6(0.06)}$&	$\bold{12(0.12)}$\\
16705  &9(1.03)&15(1.85)&54(5.65)&21(4.25)&45(8.95)&71(14.0)&$\bold{4(0.34)}$&	$\bold{6(0.45)}$&	$\bold{8(0.98)}$\\
66177  &9(4.67)&13(6.81)&49(23.9)&23(20.3)&47(41.1)&75(64.2)&$\bold{4(2.54)}$&	$\bold{6(2.04)}$&	$\bold{8(2.57)}$\\
263425&9(20.9)&13(29.4)&37(79.7)&25(107)&53(222)&77(318)&$\bold{4(6.25)}$&	$\bold{6(8.85)}$&	$\bold{8(10.9)}$\\
\hline
\end{tabular}
\end{table}
%table5
\begin{table}[!t]
\renewcommand\arraystretch{1.2}
\scriptsize
\centering
\caption{\protect\\Iterative results for Test problem 2 with $\Gamma=\Gamma_{2}$}
\begin{tabular}{c|ccc|ccc|ccc}
\hline
Method          & \multicolumn{3}{c|}{MINRES($\mathcal{P}_{D,1}$) \cite{trh1}} &  \multicolumn{3}{c|}{MINRES($\mathcal{P}_{D,2}$) \cite{jpm}} &  \multicolumn{3}{c}{$\bold{GMRES(\mathcal{P}_{T}}$)}     \\
\hline
DoF/$\beta$& $10^{-2}$  & $10^{-4}$               & $10^{-6}$  & $10^{-2}$  & $10^{-4}$               & $10^{-6}$  & $10^{-2}$  & $10^{-4}$& $10^{-6}$\\
\hline
4289    &11(0.21)&23(0.48)&95(1.87)&27(1.42)&51(2.59)&67(3.39)&$\bold{5(0.10)}$&	$\bold{8(0.15)}$&	$\bold{16(0.21)}$\\
16769  &11(1.28)&21(2.28)&89(8.93)&29(7.34)&63(15.4)&85(20.7)&$\bold{5(0.39)}$&	$\bold{8(0.64)}$&	$\bold{16(1.06)}$\\
66305  &9(5.23)&23(11.5)&83(39.0)&29(28.6)&67(64.4)&95(90.6)&$\bold{5(1.73)}$&	$\bold{8(2.54)}$&	$\bold{16(4.85)}$\\
263681&9(21.2)&23(49.9)&77(163)&31(154)&75(366)&113(554)&$\bold{5(7.34)}$&	$\bold{8(11.1)}$&	$\bold{16(22.5)}$\\
\hline
\end{tabular}
\end{table}
%table6
\begin{table}[!t]
\renewcommand\arraystretch{1.2}
\scriptsize
\centering
\caption{\protect\\Iterative results for Test problem 2 with $\Gamma=\Gamma_{3}$}
\begin{tabular}{c|ccc|ccc|ccc}
\hline
Method          & \multicolumn{3}{c|}{MINRES($\mathcal{P}_{D,1}$) \cite{trh1}} &  \multicolumn{3}{c|}{MINRES($\mathcal{P}_{D,2}$) \cite{jpm}} &  \multicolumn{3}{c}{$\bold{GMRES(\mathcal{P}_{T}}$)}     \\
\hline
DoF/$\beta$& $10^{-2}$  & $10^{-4}$               & $10^{-6}$  & $10^{-2}$  & $10^{-4}$               & $10^{-6}$  & $10^{-2}$  & $10^{-4}$& $10^{-6}$\\
\hline
4321    &15(0.31)&50(0.98)&262(4.92)&45(2.71)&81(5.00)&105(5.96)&$\bold{6(0.06)}$&	$\bold{13(0.12)}$&	$\bold{34(0.45)}$\\
16833  &17(1.84)&48(4.96)&293(29.3)&45(12.7)&99(29.0)&153(43.1)&$\bold{6(0.32)}$&	$\bold{14(0.95)}$&	$\bold{36(2.51)}$\\
66433  &13(6.43)&51(25.2)&298(144)&49(53.6)&117(131)&211(238)&$\bold{6(1.96)}$&	$\bold{14(4.39)}$&	$\bold{36(12.2)}$\\
263937&13(28.3)&43(93.0)&302(636)&51(345)&143(954)&284(1914)&$\bold{6(2.26)}$&	$\bold{14(18.7)}$&	$\bold{36(48.5)}$\\
\hline
\end{tabular}
\end{table}
The comparisons of different preconditioners are presented in Table 4-6. It is shown the preconditioner performance here is similar to that presented in Test problem 1. The cost of the GMRES($\mathcal{P}_{T}$) method is still the least and the iterations required keeps mesh size independence.

\section{Conclusion}
In this paper, a block triangular preconditioning method has been presented for the saddle point problem arising from the elliptic boundary optimal control problem with mixed boundary conditions. In the proposed method, first the saddle point problem is permuted to obtain its equivalent form. Then a preconditioner is constructed in block triangular form based on an approximation of the Schur complement. We have analyzed the eigenvalue properties of the preconditioned matrix. The effectiveness of the proposed preconditioning method was tested by comparison with other existing methods on several mixed boundary control problems.

\end{document}